\documentclass[12pt, a4paper]{article}
\usepackage[latin1]{inputenc}
\usepackage[english]{babel}
\usepackage{indentfirst}
\usepackage{amsxtra}
\usepackage{amsmath}
\usepackage{amsthm}
\usepackage{amstext}
\usepackage{amsfonts}
\usepackage{textcomp}
\usepackage{amssymb}
\usepackage{amscd}
\usepackage{color}
\usepackage[all]{xy}
\usepackage[dvips]{graphicx}
\usepackage{setspace}
\usepackage[all]{xy}

\newcommand {\F}{\mathbb{F}}

\newcommand {\la}{\lambda}
\newcommand {\al}{\alpha}

\newcommand {\Cpr}{C_{pr}^*(G)}
\newcommand {\Cprel}[1]{C_{pr}^*(G,#1)}

\newcommand{\funcao}[5]{\begin{array}{lrcl}
#1:&\!\!\!#2 & \rightarrow & #3 \\
  &\!\!\! #4 & \mapsto & #5
\end{array}}

\newtheorem{theorem}{Theorem}[section]
\newtheorem{lemma}[theorem]{Lemma}

\newtheorem{definition}[theorem]{Definition}
\newtheorem{proposition}[theorem]{Proposition}

\newtheorem{remark}[theorem]{Remark}

\begin{document}
\doublespace

\title{Partial group algebra with projections and relations}
\author{Danilo Royer}
\maketitle
\begin{abstract}We introduce the notion of {\it partial group algebra with projections and relations} and show that this C*-algebra is a partial crossed product. Examples of partial group algebras with projections and relations are the Cuntz-Krieger algebras and the unitization of  C*-algebras of directed graphs.
\end{abstract}
\maketitle




\section{Introduction}
In \cite{exelquigglaca} it was introduced the notion of {\it partial group algebra with relations}. This algebra is an universal C*-algebra where the generators are the elements of a group and the relations are, among other, the relations of a partial representation of the group. It was shown in the same paper that these algebras are examples of partial crossed products. 

In this paper we introduce the {\it partial group algebras with relations and projections}. We define this algebra in terms of generators and relations, where the generators are a family of projections and the elements of a group. We show that the algebras so defined are also partial crossed products. 

We show that the unitization of graph C*-algebras are examples of partial group algebras with relations and projections. Another class of examples of partial group algebras with projections and relations are the C*-algebras (denoted by $O_{A,B}$) which we introduce in subsection 5.1. These algebras are obtained by considering two $n\times n$ matrices $A,B$ with entries in $\{0,1\}$, and are such that if one of the matrices $A$ or $B$ is a permutation matrix then the algebra $O_{A,B}$ is the Cuntz-Krieger algebra of one matrix.

\section{Partial group algebra with projections}\label{gpap}

In this section we define the {\it partial group algebra with projections}, which is an universal C*-algebra, and we show that this algebra is a partial crossed product of a group by a commutative C*-algebra.

Let us first define the {\it partial group algebra with relations and projections}.

\begin{definition}\label{groupalgebra}Let $G$ be a group and $J$ be a set. The partial group algebra with projections, denoted $\Cpr$, is the universal C*-algebra generated by the set $\{[g]:g\in G\}$ and a family of projections $\{P_i\}_{i\in J}$ with the relations:
\begin{itemize}
\item $[e]=1$;
\item $[r^{-1}]=[r]^*$ for each $r\in G$;
\item $[r][s][s^{-1}]=[rs][s^{-1}]$ for each $r,s\in G$;
\item $P_i$ is a projection, that is, $P_i^*=P_i=P_i^2$ for each $i\in J$;
\item $[r][r]^*$ and $[s]P_j[s]^*$ commute, for all $r,s\in G$ and $i\in J$;
\item $[r]P_i[r]^*$ and $[s]P_j[s]^*$ commute, for all $r,s\in G$ and $i,j\in J$;

\end{itemize}
\end{definition}

The first, second and third relations mean that the map $G\ni g \mapsto [g]\in \Cpr$ is a partial representation of $G$. So, following 
\cite{exelpartrep}, $[r][r]^*$ and $[s][s]^*$ commute for all $r,s\in G$. Also, considering $r=s=e$, by the last relation the projections $P_i$ and $P_j$ commute, for all $i,j\in J$. In $\Cpr$ the relation $[r]^*[r][s]=[r]^*[rs]$ (similar to the third relation) also holds for each $r,s\in G$, since $[r]^*[r][s]=([s^{-1}][r^{-1}][r])^*=([s^{-1}r^{-1}][r])^*=([rs]^*[r])^*=[r]^*[rs]$.

\section{The partial group algebra as a partial crossed product}

Another way to look to the partial group algebra with projections is via the partial crossed product theory  (in the sense of \cite{exelpartrep}, \cite{mclanahan}). From now on to the end of this section, we will concentrate our attention to show that the algebra $\Cpr$ is a partial crossed product of $G$ by a commutative C*-algebra.

Let us begin by defining the topological space which gives rise to the desired commutative C*-algebra.

First add a new element to $J$, called $0$, that is, consider the new set $I$ as the disjoint union $I=J\bigsqcup \{0\}$. We will use the symbol $\bigsqcup$ for disjoint unions.

For each $i\in I$ consider a distinct copy $G^i$ of $G$ and consider the disjoint union $\bigsqcup\limits_{i\in I} G^i$. 
Given $i\in I$ with $i\neq 0$ (that is, $i\in J$) we denote the elements of $G^i$ by $g^i$. The set $G^0$ will be identified with $G$ and the elements of $G^0$ will be denoted by $g$.
Consider the topological space $$\{0,1\}^{\bigsqcup\limits_{i\in I}G^i}$$ with the product topology. This is a compact Hausdorff space. The elements of this space may be identified with the subsets of $\bigsqcup\limits_{i\in I}G^i$ via the identification $$\{0,1\}^{\bigsqcup\limits_{i\in I} G^i}\ni \widetilde{\xi}=(\widetilde{\xi}_x)_{x\in\bigsqcup\limits_{i\in I} G^i}\longleftrightarrow \xi\in\bigsqcup\limits_{i\in I} G^i$$
where $\xi=\{x\in \bigsqcup\limits_{i\in I} G^i: \widetilde{\xi}_x=1\}$.

Consider the set

$$X_G=\left\{\xi\in\{0,1\}^{\bigsqcup\limits_{i\in I}G^i}: 
\begin{array}{l}- e\in \xi\\- \text{ if } r^i\in \xi \text{ then } r\in \xi
\end{array}\right\}.$$

Note that $X_G$ is closed in $\{0,1\}^{\bigsqcup\limits_{i\in I}G_i}$ and so $X_G$ is compact. 


The C*-algebra $C(X_G)$ is the commutative $C^*$-algebra which will be used to define a partial action to obtain the partial crossed product. The partial action will be introduced later.

To show that the algebra $\Cpr$ is a partial crossed product of $G$ by a commutative C*-algebra, it is important to identify a commutative C*-algebra in $\Cpr$. Note that the subalgebra of $\Cpr$ generated by $[r][r]^*$ and $[r]P_i[r]^*$, with $r\in G$ and $i\in J$ is commutative, by definition.

With this subalgebra in mind, let us consider the  universal C*-algebra $B$ generated by elements $\{Q_{r^i}\}_{r^i\in\bigsqcup\limits_{i\in I}G^i}$ with the relations:
\begin{itemize}
\item $Q_{r^i}$ are commutative projections;
\item $Q_e=1$, where $e\in G$ is the neutral element;
\item $Q_{r^i}Q_r=Q_{r^i}$ for each $r\in G$ and $i\in I$.
\end{itemize} 
Note that $B$ is commutative and let $\widehat{B}$ be the spectrum of $B$, which is compact since $B$ is unital. 

\begin{lemma}
$\widehat{B}$ and $X_G$ are homeomorphic.
\end{lemma}

\proof Given an element $\xi\in X_G$, define $T_\xi$ on the generators of $B$ by $T_\xi(Q_{r^i})=[r^i\in \xi]$ (the notation $[r^i\in \xi]$ means: $[r^i\in \xi]=1$ if $r^i\in \xi$ and $[r^i\in \xi]=0$ if $r^i\notin \xi$). By the universal property of $B$, $T_\xi$ may be extended to $B$, and so $T_\xi \in \widehat{B}$. So we obtain an injective function
$$\funcao{T}{X_G}{\widehat{B}}{\xi}{T_\xi}.$$ For each $\psi\in \widehat{B}$  define $$\nu_\psi=\{r^i\in \bigsqcup\limits_{i\in I}G^i\,\,:\,\,\psi(Q_{r^i})=1\}.$$ Note that $e\in \nu_\psi$ because $\psi(Q_e)=\psi(1)=1$, and if $r^i\in \nu_\psi$ then $\psi(Q_{r^i})=1$, and since $\psi(Q_{r^i})=\psi(Q_{r^i})\psi(Q_r)$, it remains that $\psi(Q_r)=1$, from where $r\in \nu_\psi$. So $\nu_\psi \in X_G$ and $T_{\nu_\psi}=\psi$, and this shows that $T$ is surjective. To see that $T$ is continuous, let $(\xi_m)_{m\in M}\subseteq X_G$ be a net such that  $\xi_m\rightarrow \xi\in X_G$ ($M$ is a directed set). Then, given $r^i\in \bigsqcup\limits_{i\in I}G^i$, there exists $m_0\in M$ such that for each $m\geq m_0$, $r^i\in \xi_m$ if and only if $r^i\in \xi$. Therefore, for each $m\geq m_0$, 

$$T_{\xi_m}(Q_{r^i})=[r^i\in \xi_m]=[r^i\in \xi]=T_{\xi}(Q_{r^i}),$$

and this shows that $T_{\xi_m}(Q_{r^i})\rightarrow T_{\xi}(Q_{r^i})$. Then $T_{\xi_m}(d)\rightarrow T_{\xi}(d)$ for $d$ in a dense set in $B$. Since $T_{\xi_m}$ are uniformly limited it follows that $T_{\xi_m}\rightarrow T_\xi$ in $\widehat{B}$. So, $T$ is bijective and continuous. Since $\widehat{B}$ and $X_G$ are compact and Hausdordff then $T$ is a homeomorphism.
\endproof

By using the homeomorphism introduced in the previous lemma, we obtain the *-isomorphism $\Phi_0:C(\widehat{B})\rightarrow C(X_G)$ defined by $\Phi_0(f)=f\circ T$ for each $f\in C(\widehat{B})$. Let $\Gamma:B\rightarrow C(\widehat{B})$ be the Gelfand *-isomorphism, that is, $\Gamma(b)(\psi)=\psi(b)$ for each $b\in B$ and $\psi\in \widehat{B}$. Now, let $\Phi:B\rightarrow C(X_G)$ be the *-isomorphism $\Phi_0\circ \Gamma$, and note that 
$$\Phi(Q_{r^i})(\xi)=(\Phi_0(\Gamma(Q_{r^i}))(\xi)=\Gamma(Q_{r^i})(T(\xi))=T_\xi(Q_{r^i})=[r^i\in \xi],$$ where $[r^i\in \xi]=1$ if $r^i\in \xi$ and $[r^i\in \xi]=0$ if $r^i\notin \xi$. Denoting by  $1_{r_i}$ the characteristic function $1_{r^i}(\xi)=[r^i\in \xi]$, we get $\Phi(Q_{r^i})=1_{r^i}$.

Let $\varphi_0: C(X_G)\rightarrow B$ be the inverse *-isomorphism of $\Phi$. By the universal property of $B$ there exists a *-homomorphism $\varphi_1:B\rightarrow \Cpr$ such that $\varphi_1(Q_{r^i})=[r]P_i[r]^*$ for each $r\in G$ and $i\in J$ and $\varphi_1(Q_{r})=[r][r]^*$ for each $r\in G$. Let
 
\begin{equation}\label{homom}
\varphi:C(X_G)\rightarrow \Cpr 
\end{equation}

be the composition $\varphi_1\circ\varphi_0$. Note that $\varphi(1_{r^i})=[r]P_i[r]^*$ for each $r\in G$ and $i\in J$ and $\varphi(1_{r})=[r][r]^*$ for each $r\in G$. 

Let us now introduce the necessary partial action of $G$ in $C(X_G)$ to obtain a partial crossed product.

For each $t\in G$, consider $X_t\subseteq X_G$, 
$$X_t=\{\xi \in X_G\,\,:\,\,t\in \xi\}$$ which is clopen in $X_G$. Given an element $\xi \in X_{t^{-1}}$ we define the set $t\xi$ by $t\xi=\{(tr)^i\,\,:\,\,r^i\in \xi\}$. Note that in this case $t\xi$ is an element of $X_t$. Then we define the map  

$$\funcao{\theta_t}{X_{t^{-1}}}{X_t}{\xi}{t\xi}$$
which is a homeomorphism.

For each $t\in G$ define the *-isomorphism 

$$\funcao{\al_t}{C(X_{t^{-1}})}{C(X_t)}{f}{f\circ \theta_{t^{-1}}}.$$ Then $(\{C(X_t)\}_{t\in G}, \{\al_t\}_{t\in G})$ is a partial action of $G$ on $C(X_G)$ (in the sense of \cite{exelpartrep}, \cite{mclanahan}).

In the following lemma, which is a technical but a useful result, $\prod1_{r^i}$ means a product of finitely many characteristic functions, where $r^i\in \bigsqcup\limits_{i\in I}G^i$. These characteristic functions are elements of $C(X_G)$. Recall that since $r^0$ is identified with $r$, in the product $\prod1_{r^i}$, characteristic functions like $1_r$, with $r\in G$ are also allowed. 
\begin{lemma}
\begin{enumerate}\label{l1}
\item $span\{1_t\prod1_{r^i}\}$ is a dense *-algebra in $C(X_t)$.
\item $\al_t(1_{t^{-1}}1_{r^i})=1_t1_{(tr)^i}$.
\end{enumerate}
\end{lemma}

\proof The first item follows from the Stone-Weierstrass theorem. For the second one, consider $\xi \in X_G$. If $\xi \in X_t$ then $$\al_t(1_{t^{-1}}1_{r^i})(\xi)=(1_{t^{-1}}1_{r^i})(t^{-1}\xi)=1_{t^{-1}}(t^{-1}\xi)1_{r^i}(t^{-1}\xi)=[r^i\in t^{-1}\xi]=$$
$$=[(tr)^i\in \xi]=1_{(tr)^i}(\xi)=1_t(\xi)1_{(tr)^i}(\xi)=(1_t1_{(tr)^i})(\xi).$$ If $\xi \notin X_t$ then 
$\al_t(1_{t^{-1}}1_{r^i})(\xi)=0=(1_t1_{(tr)^i})(\xi)$.
\endproof


Define the map $\pi:G\rightarrow \Cpr$ by $\pi(g)=[g]$ and note that (by the definition of $\Cpr$) $\pi$ is a partial representation of $G$.


\begin{proposition}\label{covrep}
The pair $(\varphi, \pi)$ is a covariant representation of $(C(X_G),\alpha, G)$, that is, $\pi$ is a partial representation of $G$, $\varphi$ is a *-homomorphism of $C(X_G)$, and  $\varphi(\al_t(a))=\pi(t)\varphi(a)\pi(t)^*$ for each $a\in C(X_t)$ and $t\in G$.
\end{proposition}

\proof All we need to do is to show that $\varphi(\al_t(a))=\pi(t)\varphi(a)\pi(t)^*$ for each $a\in C(X_t)$ and $t\in G$.

Let $t\in G$. Then, given $r^i\in \bigsqcup\limits_{i\in I}G^i$,   
$$\varphi(\al_{t}(1_{t^{-1}}1_{r^i}))=\varphi(1_t1_{(tr)^i})=[t][t]^*[tr]P_i[tr]^*=[t][t]^*[tr]P_i[tr]^*[t][t]^*=$$
$$=[t][r]P_i[r]^*[t]^*=[t][t]^*[t][r]P_i[r]^*[t]^*=\pi(t)\varphi(1_{t^{-1}}1_{r^i})\pi(t)^*.$$

By Lemma $\ref{l1}$, $\varphi(\al_t(b))=\pi(t)\varphi(b)\pi(t)^*$ for each $b$ in a dense subset of $C(X_t)$. The result follows by continuity of $\varphi$ and $\alpha_t$.
\endproof



We recall that the partial crossed product of $G$ and $C(X_G)$ by the partial action $\alpha$, (see \cite{exelcircleactions}, \cite{mclanahan}), denoted by $C(X_G)\rtimes_\al G$,  is the enveloping C*-algebra of the *-algebra $\left\{\sum\limits_{g}a_g\delta_g\,:\, a_g\in C(X_g)\right\}$ endowed with: the pointwise sum, the product 
$(\sum\limits_ga_g\delta_g)(\sum\limits_h b_h\delta_h)=\sum\limits_{h,g}\alpha_g(\alpha_{g^{-1}}(a_g)b_{g^{-1}h})\delta_h$, the involution $(\sum\limits_ga_g\delta_g)^*=\sum\limits_g\alpha_{g^{-1}}(a_g^*)\delta_{g^{-1}}$ and $\|\sum\limits_{g}a_g\delta_g\|_1=\sum\limits_{g}\|a_g\|$.

Now we are ready to prove that the partial group algebra with projections  $\Cpr$ introduced at the beginning of this section is isomorphic to the partial crossed product $C(X_G)\rtimes_\al G$.

\begin{theorem}\label{t1}
There exists a *-isomorphism $\phi:\Cpr\rightarrow C(X_G)\rtimes_\al G$ such that $\phi([r][r]^*)=1_r\delta_e$ and $\phi(P_i)=1_{e^i}\delta_e$.
\end{theorem}

\proof We will prove that there exist two *-homomorphism, one from $\Cpr$ to $C(X_G)\rtimes_\al G$, and other from $C(X_G)\rtimes_\al G$ to $\Cpr$, and that these *-homomorphisms are mutually inverse.
First, define $$\varphi\times \pi:\left\{\sum a_g\delta_g\,:\, a_g\in C(X_g)\right\}\rightarrow \Cpr$$ by $(\varphi\times\pi)(\sum a_g\delta_g)=\sum \varphi(a_g)\pi(g)$ where $\varphi$ and $\pi$ are the maps of Proposition \ref{covrep}.
It is easy to show, by using Proposition \ref{covrep} and direct calculations, that $\varphi\times \pi$ is linear, multiplicative, self-adjoint and contractive, and so $\varphi\times \pi$ extends to $C(X_G)\rtimes_\alpha G$.

On the other hand, by the universal property of $\Cpr$ there exists a *-homomorphism $\phi:\Cpr\rightarrow C(X_G)\rtimes_\alpha G$ such that $\phi([r])=1_r\delta_r$ and $\phi(P_i)=1_{e^i}\delta_e$, for each $r\in G$ and $i\in J$. 

Note that
$$\phi([r][r]^*)=1_r\delta_r1_r^{-1}\delta_r^{-1}=1_r\delta_e.$$

To show that $\phi$ is a *-isomorphism, we show that $\varphi\times \pi$ is the inverse of $\phi$. 

For each $t\in G$ and $i\in J$,  $$(\phi\circ(\varphi\times\pi))(1_t1_{r^i}\delta_t)=\phi(\varphi(1_t1_r^i)\pi(t))=$$
$$=\phi([t][t]^*[r]P_i[r]^*[t])=1_t\delta_e1_r\delta_r1_{e^i}\delta_e1_{r^{-1}}\delta_{r^{-1}}1_t\delta_t=1_t1_r1_{r^i}\delta_t.$$
Recall that given $\xi \in X_G$, if $r^i\in \xi$ then $r\in \xi$. Therefore, $1_r1_{r^i}=1_{r^i}$ in $C(X_G)$. Then  
$1_t1_r1_{r^i}\delta_t=1_t1_{r^i}\delta_t.$ It follows from Lemma \ref{l1} that $\phi\circ(\varphi\times\pi)$ is the identity homomorphism of $C(X_G)\rtimes_\al G$.

On the other hand, by direct calculations, $((\varphi\times\pi)\circ\phi)([r])=[r]$ and $((\varphi\times\pi)\circ\phi)(P_i)=P_i$, and so  
$(\varphi\times\pi)\circ\phi$ is the identity homomorphism of $\Cpr$. 

\endproof


\section{The partial group algebra with projections and relations}

The partial group algebra with projections is, roughly speaking, a C*-algebra generated by partial isometries and projections. One of the most celebrated C*-algebra generated by partial isometries is the Cuntz-Krieger algebra $O_A$ for a $n\times n$ matrix $A$ with entries in $\{0,1\}$. If we desire (for example) to show that the Cuntz-Krieger algebras are examples of our construction, we need to identify some projections in these algebras. Let $\{S_1,...,S_n\}$ be the generators of $O_A$, and let $p_i$ be the projection $S_iS_i^*$. So, in $O_A$ the relations $p_i=S_iS_i^*$, $S_i^*S_i=\sum\limits_{j=1}^na_{i,j}p_j$ and $\sum\limits_{i=1}^np_i=1$ hold. To identify $O_A$ with some partial group algebra with projections, it seems to be reasonable to choose the group as being the free group $\F_n$ (generated by $\{1,2,...,n\}$), and $n$ projections $P_1,...,P_n$, so that the partial isometries $S_i$ and the projections $p_i$ of $O_A$ identifies with the partial isometries $[i]$ and the projections $P_i$ of $C_{pr}^*(\F_n)$. In $O_A$, the equalities $S_i^*S_i=\sum\limits_{j=1}^np_i$, $S_iS_i^*=p_i$ and $\sum\limits_{i=1}^np_i=1$ hold, but nothing guarantees that the relation $[i][i]^*=P_i$ holds in $C_{pr}^*(\F_n).$ So, it is natural to consider the {\it partial group algebra with projections and relations}. In this section, we define the {\it partial group algebra with projections and relations} and show that this algebras are also partial crossed products of a group by a commutative C*-algebra.

Let us begin by defining the {\it partial group algebra with projections and relations}. 

Let $G$ be a group, $J$ be a set and let $X_G$ be the topological space as in  Section \ref{gpap}. Consider a set $R\subseteq C(X_G)$ such that each element $f\in R$ is of the form $f=\sum\limits_j\la_j\prod\limits_{r,i}1_{r^i}$, that is, $f$ is a linear combination of products of finitely many characteristic functions $1_{r_i}$, where $r^i\in \bigsqcup\limits_{i\in I}G^i$. Recall that $I=J\bigsqcup\limits\{0\}$, and since we are identifying $G^0$ with $G$, functions like $1_r$ with $r\in G$ are also allowed in the products $\prod\limits_{r,i}1_{r^i}$. 

\begin{definition} Let $G$ be a group, $J$ be a set and let $X_G$ be the topological space as in  section \ref{gpap}. In $C(X_G)$, consider a set $R$ as above.
The partial group C*-algebra with relations $R$ and projections $\{P_i\}_{i\in J}$, which we denote $\Cprel{R}$, is the universal C*-algebra generated by $\{[g]:g\in G\}$ and $\{P_i\}_{i\in J}$ with the relations of the Definition \ref{groupalgebra} and the relations
\item $\sum\limits_j\la_j\prod\limits_{r,i}[r]P_ i[r]^*=0$ for all $\sum\limits_j\la_j\prod\limits_{r,i}1_{r^i}\in R$. Here, $[r]P_0[r]^*$ means $[r][r]^*$.

\end{definition}

By definition, the algebra $\Cprel{R}$ is a quotient of $\Cpr$. If $R=\emptyset$ then $\Cprel{R}=\Cpr$.

As in the previous section, we will show that the algebras $\Cprel{R}$ are partial crossed products of $G$ by a commutative algebra. 

First we let $R\subseteq C(X_G)$ and consider the ideal of $C(X_G)$ generated by $$\{ \al_t(1_{t^{-1}}f):f\in R, t\in G\}.$$ This ideal is of the form $C_0(V)$, for some open subset $V$ of $X_G$.

\begin{proposition}
Define $$\Omega_R=\{\xi\in X_G\,:\, f(t^{-1}\xi)=0\,\forall\, t\in \xi\cap G,\,f\in R\}.$$ Then $X_G\setminus V=\Omega_R$, (therefore $C(\Omega_R)\simeq \frac{C(X_G)}{C_0(V)}$).
\end{proposition}
 
\proof
Let $\xi\in X_G\setminus V$. Then, given $t\in \xi\cap G$ and $f\in R$, since $\al_t(1_{t^{-1}}f)\in C_0(V)$ then $0=\al_t(1_{t^{-1}}f)(\xi)=(1_{t^{-1}}f)(t^{-1}\xi)=f(t^{-1}\xi)$. Therefore $\xi\in \Omega_R$.

On the other hand, let $\xi \in \Omega_R$. Then, given $f\in R$, $\al_t(1_{t^{-1}}f)(\xi)=0$ if $t\notin \xi$, and if $t\in\xi$ then $\al_t(1_{t^{-1}}f)(\xi)=1_{t^{-1}}f(t^{-1}\xi)=f(t^{-1}\xi)=0$. So, $g(\xi)=0$ for each $g\in C_0(V)$,  and so $\xi \notin V$.  
\endproof

The set $\Omega_R$ is a compact set in $X_G$ and is invariant under the partial action defined for each $t\in G$ by $X_{t^{-1}}\ni \xi\mapsto t\xi\in X_t$. Then we may consider the restriction of the partial action  $\theta$ to $\Omega_R$, that is, we may consider the partial action $\theta_t:X_{t^{-1}}\cap \Omega_R\rightarrow X_t\cap \Omega_R$, which we also call $\theta$. Then, by defining $D_t=C(X_t\cap \Omega_R)$, and $$\funcao{\beta}{D_{t^{-1}}}{D_t}{f}{f\circ \theta_t}$$ we obtain a partial action of $G$ on $C(\Omega_R)$. 

From now on to the end of this section let us consider only subsets $R$ of $C(X_G)$ such that all the elements of $R$ are of the form 
$\sum\limits_{j=1}^n\la_j\prod\limits_{r,i}1_{r^i}$.

\begin{theorem}
There is a *-isomorphism $\phi:\Cprel{R}\rightarrow C(\Omega_R)\rtimes_\beta G$ such that $\phi([r])=1_r\delta_r$ and $\phi(P_i)=1_{e^i}\delta_e$.  
\end{theorem} 

\proof
By the universal property of $\Cprel{R}$, there is a *-homomorphism 
$\phi:\Cprel{R}\rightarrow C(\Omega_R)\rtimes_\beta G$ such that $\phi([r])=1_r\delta_r$ and $\phi(P_i)=1_{r^i}\delta_e$ for each $r\in G$ and $i\in J$. We will only verify that $\phi(f)=0$ for each $f\in R$. The verification of the other relations are left to the reader.
Let $f=\sum\limits_j\la_j\prod\limits_{r,i}1_{r^i}\in R$. Then, 
$$\phi(f)=\sum\limits_j\la_j\prod\limits_{r,i}\phi([r])\phi(P_i)\phi([r])^*=\sum\limits_j\la_j\prod\limits_{r,i}1_r\delta_r1_{e^i}\delta_e1_{r^{-1}}\delta_{r^{-1}}=$$
$$=\left(\sum\limits_j\la_j\prod\limits_{r,i}1_r1_{r^i}\right)\delta_e.$$
Since $1_r1_{r^i}=1_{r^i}$ it follows that $\sum\limits_j\la_j\prod\limits_{r,i}1_r1_{r^i}=\sum\limits_j\la_j\prod\limits_{r,i}1_{r^i}=0$ in $\Omega_R$ and so $\phi(f)=0$.

On the other hand, consider the *-homomorphism $\overline{\varphi}:C(X_G)\rightarrow \Cprel{R}$ obtained by composing the homomorphism (\ref{homom}) with the quotient map from $\Cpr$ to $\Cprel{R}$. Note that if $f=\sum\limits_j\la_j\prod\limits_{r,i}1_{r^i}\in R$ then 

$$\overline{\varphi}(\al_t(1_{t^{-1}}f))=\overline{\varphi}\left(\sum\limits_j\la_j\prod\limits_{r,i}(1_t1_{(tr)^i})\right)=$$
$$=\sum\limits_j\la_j\prod\limits_{r,i}[t][t]^*[tr]P_i[tr]^*=\sum\limits_j\la_j\prod\limits_{r,i}[t][t]^*[tr]P_i[tr]^*[t][t]^*=$$
$$=\sum\limits_j\la_j\prod\limits_{r,i}[t][r]P_i[r]^*[t]^*=[t]\left(\sum\limits_j\la_j\prod\limits_{r,i}[r]P_i[r]^*\right)[t]^*=0.$$
So $\overline{\varphi}(g)=0$ for each $g$ in the ideal $C_0(V)$. Since $C(\Omega_R)\simeq \frac{C(X_G)}{C_0(V)}$, there is a *-homomorphism $\psi:C(\Omega_R)\rightarrow \Cprel{R}$ such that $\psi(1_{r^i})=[r]P_i[r]$. 
Consider the map $\pi:G\rightarrow \Cprel{R}$ defined by $\pi(g)=[g]$, which is a partial representation. Then $(\psi, \pi)$ is covariant and the map
$$\psi\times\pi:\left\{\sum a_g\delta_g\,:\, a_g\in D_g,\,g\in G\right\}\rightarrow \Cprel{R}$$

 defined by $(\psi\times\pi)(\sum a_g\delta_g)=\sum\psi(a_g)\pi(g)$
is a contractive *-homomorphism and extends to $C(\Omega_R)\rtimes_\beta G$. The *-homomorphisms $\phi$ and $\psi\times \pi$ are inverses of each other.\endproof

\section{Examples}

\subsection{The Cuntz-Krieger C*-algebras}

Recall that the celebrated Cuntz-Krieger C*-algebra $O_D$ of a $n\times n$ matrix $D=(d_{i,j})_{i,j}$ with entries in $\{0,1\}$, is the universal C*-algebra generated by $n$ partial isometries $S_i$ with the relations $\sum\limits_{i=1}^nS_iS_i^*=1$ and $S_i^*S_i=\sum\limits_{j:d_{i,j=1}}S_jS_j^*$.

Let $A=(a_{i,j})_{i,j}$ and $B=(b_{i,j})_{i,j}$ be two matrices of order $n\times n$ such that $a_{i,j}, b_{i,j}\in \{0,1\}$. Let $\F_n$ be the free group with $n$ generators. 

Let $O_{A,B}$ be the partial group C*-algebra with relations and projections $C_{pr}^*(\F_n,R)$ generated by the group $\F_n$ and projections $P_1,...,P_n$ with the relations: 

\begin{itemize}
\item $\sum\limits_{i=1}^nP_i=1$.
\item $[i]^*[i]=\sum\limits_{j=1}^na_{i,j}P_j$, for each $i\in \{1,...,n\}$.
\item $[i][i]^*=\sum\limits_{j=1}^nb_{i,j}P_j$, for each $i\in \{1,...,n\}$.
\item $[r][r]^*[rs][rs]^*=[rs][rs]^*$ for all $r,s\in \F_n$ such that $|rs|=|r|+|s|$, (where $|r|$ means the length of $r$, that is, the quantity of generators of $r$ in its reduced form).
\end{itemize}

It is  easy to see that $O_{A,B}$ and $O_{B,A}$ are isomorphic C*-algebras, where one isomorphism is $O_{A,B}\ni x\mapsto x^*\in O_{B,A}$.

The celebrated Cuntz-Krieger algebras $O_A$ may be seen as algebras $O_{A,B}$ following the next proposition.

\begin{proposition}
If one of the matrices $A$ or $B$ is a permutation matrix then $O_{A,B}$ is a Cuntz-Krieger algebra. More precisely, if $B$ is a permutation matrix then $O_{A,B}$ is *-isomorphic to the Cuntz-Krieger algebra $O_{\widetilde{A}}$, where $\widetilde{A}=B^tA$. 
\end{proposition}

\proof First write $B=(b_{i,j})_{i,j}$ and note that for each $i$ there is a unique $j$ such that $b_{i,j}=1$ (and for each $j$ there is a unique $i$ with $b_{i,j}=1$). Moreover, note that the line $j$ of $\widetilde{A}$ is equal to the line $i$ of $A$, that is, $\widetilde{a}_{j,k}=a_{i,k}$ for each $k$, where $i,j$ are such that $b_{i,j}=1$.

Define $\pi:\{1,...,n\}\rightarrow O_{\widetilde{A}}$ by $\pi(i)=S_j$, where $i,j$ are such that $b_{i,j}=1$. Here, $S_1,...,S_n$ are the generators of $O_{\widetilde{A}}$. By \cite[5.2]{amenabilityforfell}, $\pi$ extends to a partial representation of $\F_n$, which we also call $\pi$. Define $\varphi:\F_n\cup\{P_1,...,P_n\}\rightarrow O_{\widetilde{A}}$ by $\varphi(r)=\pi(r)$ for $r\in \F_n$ and $\varphi(P_i)=S_iS_i^*$. 
To extend $\varphi$ to $O_{A,B}$ we need to verify the relations defining $O_{A,B}$. First note that the relation $\varphi(rs)=\varphi(r)\varphi(s)$ for $r,s\in \F_n$ with $|rs|=|r|+|s|$ follows by the definition of $\varphi$. The equality $\sum\limits_{i=1}^n\varphi(P_i)=1$ is trivial. Let $i\in \{1,...,n\}$ and choose $j$ such that $b_{i,j}=1$ (such a $j$ is unique). Then 
$$\varphi(i)^*\varphi(i)=S_j^*S_j=\sum\limits_{k=1}^n\widetilde{a}_{j,k}S_kS_k^*=\sum\limits_{k=1}^na_{i,k}S_kS_k^*=\sum\limits_{k=1}^na_{i,k}\varphi(P_k),$$
and 
$$\varphi(i)\varphi(i)^*=S_jS_j^*=\varphi(P_j)=\sum\limits_{k=1}^nb_{i,k}\varphi(P_k).$$

So, we obtain a *-homomorphism (which we also call $\varphi$) $\varphi:O_{A,B}\rightarrow O_{\widetilde{A}}$.

To obtain the inverse of $\varphi$ we first define $\psi:\{S_1,...,S_n\}\rightarrow O_{A,B}$ by $\psi(S_j)=[i]$, where $i,j$ are such that $b_{i,j}=1$, and then extend it to $O_{\widetilde{A}}$ by using the universal property of $O_{\widetilde{A}}$.

\endproof









\subsection{The graph algebra}

Consider a directed graph $E=(E^0, E^1, r,s)$ and the graph C*-algebra  $C^*(E)$. Following the convention of \cite{raeburn}, $C^*(E)$ is the universal C*-algebra generated by mutually orthogonal projections $\{Q_v\}_{v\in E^0}$ and by partial isometries $\{S_a\}_{a\in E_1}$ with the relations:

\begin{itemize}
\item $S_a^*S_b=0$ for each $a,b\in E^1, a\neq b$, that is, the partial isometries have orthogonal ranges,  
\item $S_a^*S_a=Q_{s(a)}$ for each $a\in E^1$,
\item $S_aS_a^*\leq Q_{r(a)}$ for each $a\in E^1$,
\item $Q_{v}=\sum\limits_{b\in r^{-1}(v)}S_bS_b^*$ if $0<\#\{r^{-1}(v)\}<\infty $.
\end{itemize}

Since the graph algebras are not always unital and the partial group algebras are unital $C^*$-algebras, it is obviously not reasonable to expect that graph algebras are *-isomorphic to partial group algebras with projections and relations. However, we show that $C^*(E)$ is *-isomorphic to a subalgebra of a  group partial algebra with projections and relations. 

Consider the free group $\F$ generated by $E^1$. Given an element $0\neq t\in \F$, write $t=t_1...t_n$ in his reduced form and define the element $\pi(t_1)...\pi(t_n)\in C^*(E)$ where $\pi(t_i)=S_a$ if $t_i=a\in E^1$ and $\pi(t_i)=S_a^*$ if $t_i=a^{-1}\in (E^1)^{-1}$. Define the map $\pi:\F\rightarrow\widetilde{C^*(E)}$ by $\pi(t_1...t_n)=\pi(t_1)...\pi(t_n)$ and $\pi(e)=1$. 

Let us show, in the next pages, that the map $\pi$ is a partial group representation. 

Let $I\subseteq C^*(E)$ be the semigroup generated by $S_a, S_b^*:a,b\in E^1$. 

Note that for each $c\in I$, $cc^*\leq Q_u$, for some $u\in E^0$. In fact, if $c=S_a\tilde{c}$ then $cc^*\leq S_aS_a^*\leq Q_{r(a)}$, and if $c=S_a^*\tilde{c}$ then $cc^*\leq S_a^*S_a = Q_{s(a)}$.

\begin{proposition}\label{partialisometry}
For each $c\in I$, $c$ is a partial isometry. Moreover, for each $c,d\in I$, $cc^*$ commutes with $dd^*$.
\end{proposition}

\proof
First we prove the following claim.

{\it Claim: If $c,d\in I$ are such that $|c|=1$ then $cc^*$ commutes with $dd^*$.}

Let $c=S_a$. If $d=S_br$ then $cc^*dd^*=S_aS_a^*S_brr^*S_b^*=[a=b]S_arr^*S_a^*=S_brr^*S_b^*S_aS_a^*=dd^*cc^*$. If $d=S_b^*r$ then $d=S_{d_n}^*...S_{d_1}^*$ with $d_i\in E^1$, and so $dd^*=Q_{s(d_n)}$ and since $cc^*\leq Q_{r(a)}$ then $cc^*dd^*=cc^*Q_{r(a)}Q_{s(d_n)}=[r(a)=s(d_n)]cc^*=[r(a)=s(d_n)]Q_{s(d_n)}Q_{r(a)}cc^*=dd^*cc^*$.
  
Let $c=S_a^*$. Then $cc^*=S_a^*S_a=Q_{s(a)}$ and in this case $cc^*$ and $dd^*$ commute for each $d\in I$ because $dd^*\leq Q_u$ for some $u\in E^0$.
 
This proves the claim. 

We will prove the proposition by using induction over $|c|$. 
Suppose $c$ is a partial isometry for each $c\in I$ with $|c|<n$. Let $d\in I$ with $|d|=n$. Write $d=tc$, where $t,c\in I$, $|t|=1$ and $|c|=n-1$. By the Claim, $t^*t$ commutes with $cc^*$ and $c$ is a partial isometry. Then 
$$dd^*d=tcc^*t^*tc=tt^*tcc^*c=tc=d,$$ and so $d$ is a partial isometry.

It follows that each element of $I$ is a partial isometry. 

By \cite[5.4]{exelpartrepamenable}, for each $c,d\in I$, $cc^*$ and $dd^*$ commute.




\endproof

Let $\widetilde{C^*(E)}$ denote the unitization of $C^*(E)$.

By the previous proposition and by \cite[5.4]{exelpartrepamenable} it follows that the map $\pi:\F\rightarrow \widetilde{C^*(E)}$ defined by $\pi(e)=1_{\widetilde{C^*(E)}}$ and $\pi(t_1...t_n)=\pi(t_1)...\pi(t_n)$ (where $t_1...t_n$ is in the reduced form) is a partial group representation.

Consider the topological space $X_{\F}$ as in section \ref{gpap}. Fix the subset $R_E$ of $C(X_\F)$ (the set of the relations) as being the set of the following functions:
\begin{itemize}
\item $1_{e^u}1_{e^v}$ : $u,v\in E^0$, with $u\neq v$ ($e$ is the neutral element of $\F$);
\item $1_{a^{-1}}-1_{e^{s(a)}}$ : $a\in E^1$;
\item $1_a1_b$ : $a,b\in E^1$ with $a\neq b$;
\item $1_{e^v}-\sum\limits_{a\in r^{-1}(v)}1_a$\,\,\, for $v\in E^0$ with\,\, $0<\#r^{-1}(v)<\infty$;
\item $1_{rs}1_r-1_{rs}$ \,\,\, for $r,s\in \F$ with $|rs|=|r|+|s|$;
\item $1_a1_{e^{r(a)}}-1_a$ : $a\in E^1$.
\end{itemize}

\begin{theorem}
The partial group C*-algebra with projections and relations $C^*(\F, R_E)$ generated by $\{[s]:s\in \F\}$ and $\{P_i\}_{i\in E^0\cup \{0\}}$ with relations $\sum\prod[r]P_i[r]^*=0$ for all $\sum\prod 1_{r^i}\in R_E$ is *-isomorphic to $\widetilde{C^*(E)}$.
\end{theorem}

\proof First note that $C^*(\F,R_E)$ is the universal $C^*$-algebra generated by $\{[a]:a\in \F\}$ and $\{P_i\}_{i\in E^0\cup\{0\}}$ with the following relations: 

\begin{itemize}
\item $[e]=1$, $[a^{-1}]=[a]^*$, and $[r][s][s^{-1}]=[rs][s^{-1}]$ for all $r,s\in \F$;
\item $[r]P_i[r]^*$ and $[s]P_j[s]^*$ commute, for all $r,s\in \F$ and $i,j\in E^0\cup\{0\}$;
\item $P_0=1$ and  $\{P_u\}_{u\in E^0}$ are orthogonal projections; 
\item $[a]^*[a]=P_{s(a)}$ for all $a\in E^1$; 
\item $\{[a][a]^*\}_{a\in E^1}$ are mutually orthogonal; 
\item $P_v=\sum\limits_{a\in r^{-1}(v)}[a][a]^*$ if $0<\#r^{-1}(v)<\infty$; 
\item $[rs][rs]^*[r][r]^*=[rs][rs]^*$ if $r,s\in \F$ with $|rs|=|r|+|s|$; 
\item $[a][a]^*P_{r(a)}=[a][a]^*$ for all $a\in E^1$.
\end{itemize}
By the universal property of $C^*(E)$ there is a *-homomorphism $\psi_0:C^*(E)\rightarrow C^*(\F, R_E)$ such that $\psi_0(S_a)=[a]$ and $\psi_0(Q_v)=P_v$. 
So, there is a *-homomorphism $$\funcao{\psi}{\widetilde{C^*(E)}}{C^*(\F,R_E)}{x+\alpha 1_{\widetilde{C^*(E)}}}{\psi_0(x)+\alpha 1_{C^*(\F,R_E)}}.$$ 

On the other hand define $$\phi:\{[a]:a\in \F\}\cup \{P_i\}_{i\in E^0\cup \{0\}}\rightarrow \widetilde{C^*(E)}$$ by $\phi([r])=\pi(r)$ if $r\in \F$, $\phi(P_v)=Q_v$ if $v\in E^0$ and $\phi(P_0)=1_{\widetilde{C^*(E)}}$. 
We need to verify that $\phi$ satisfies the relations which define $C^*(\F, R_E)$.
By the previous proposition, $\phi_{|_{\{[r]:r\in\F\}}}$ satisfies the relations which define a partial group representation.
To see that $\phi([c])\phi(P_u)\phi([c])^*$ commutes with $\phi([d])\phi(P_v)\phi([d])^*$ (where $u,v\in E^0$ and $c,d\in \F$), it is enough to see that $\phi([c])\phi(P_v)\phi([c])^*=\delta\phi([c])\phi([c])^*$ for each $c\in \F$, where $\delta\in\{0,1\}$, because $\phi([c])\phi([c])^*$ commutes with $\phi([d])\phi([d])^*$, by [\ref{partialisometry}]. So, given $c\in \F$ ($c$ in its reduced form), write $c=\tilde{c}t$, with $|t|=1$. If $\phi([t])=S_a$ then 
$$\phi([t])\phi(P_v)=S_aQ_v=S_aS_a^*S_aQ_v=$$ $$=[s(a)=v]S_aS_a^*S_a=[s(a)=v]S_a=[s(a)=v]\phi([t]).$$
If $\phi([t])=S_b^*$ then 
$$\phi([t])\phi(P_v)=S_b^*Q_v=S_b^*S_bS_b^*Q_v=S_b^*S_bS_b^*Q_{r(b)}Q_v=$$ $$=[r(b)=v]S_b^*S_bS_b^*Q_{r(b)}=[r(b)=v]S_b^*S_bS_b^*=[r(b)=v]S_b^*=[r(b)=v]\phi([t]).$$

So, if $|t|=1$ then $\phi([t])\phi(P_v)\phi([t])^*=\delta\phi([t])\phi([t])^*$ for some $\delta\in\{0,1\}$.

Then $\phi([c])\phi(P_v)\phi([c])^*=\phi([\tilde{c}])\phi([t])\phi(P_v)\phi([t])^*\phi([\tilde{c}])^*=\delta\phi([\tilde{c}])\phi([t])\phi([t])^*\phi([\tilde{c}])^*=\delta\phi([c])\phi([c])^*$.

It is easy to show that $\phi$ satisfies the other properties which define $C^*(\F, R_E)$. So $\phi$ extends to $C^*(\F, R_E)$. The *-homomorphisms $\phi$ and $\psi$ are inverse each of other.
\endproof

\begin{remark}
Following the previous theorem, the graph algebra $C^*(E)$ is *-isomorphic to the subalgebra of $C^*(\F,R_E)$ generated by $\{P_v:v\in E^0\}$ and $\{[r]:r\in \F\setminus\{e\}\}$.
\end{remark}

\bibliographystyle{amsplain}

\end{document}